\documentclass{proc-l}

\newtheorem{theorem}{Theorem}[section]
\newtheorem{lemma}[theorem]{Lemma}

\theoremstyle{definition}
\newtheorem{definition}[theorem]{Definition}

\usepackage{pstricks}
\usepackage{graphicx}
\usepackage{amsthm}
\usepackage{amssymb}

\theoremstyle{remark}

\numberwithin{equation}{section}

\begin{document}

\title{Minimal Bounds and Members of Effectively Closed Sets}


\author{Ahmet \c{C}ev\.{I}k}
\address{Gendarmerie and Coast Guard Academy, Ankara, Turkey}
\curraddr{}
\email{a.cevik@hotmail.com}
\thanks{The author was supported during this research by the Scientific and Technological Research Council of Turkey (T\"{U}B\.{I}TAK) 1059B191500188. Special thanks to Theodore A. Slaman for many helpful discussions.}


\subjclass[2010]{Primary }

\date{}

\dedicatory{}

\commby{}


\begin{abstract}
We show that there exists a non-empty special $\Pi^0_1$ class in which no member is a minimal cover for any set, hence prove that degrees of minimal covers cannot be a basis for $\Pi^0_1$ classes.
\end{abstract}

\maketitle


%

Study of effectively closed sets, or namely $\Pi^0_1$ classes, has been a long-standing theme in classical recursion theory. Particularly the problem of determining degree theoretic complexity of members of $\Pi^0_1$ classes, going back to Kleene \cite{SK}, has resulted in a well developed theory. By compactness of the Cantor space, degree theoretic complexity of members of $\Pi^0_1$ classes also determines reals that can be defined by compactness rather than using replacement. Jockush and Soare \cite{JS1} \cite{JS2}, in their leading papers, showed some very interesting degree theoretic properties of members of $\Pi^0_1$ classes. Many of these results have come to known as {\em basis theorems} for $\Pi^0_1$ classes. A typical basis theorem states that every $\Pi^0_1$ class has a member, or a member of degree, of a particular kind. It may be the case that not every $\Pi^0_1$ class has members with the desired property. Similarly, it may be that there is a $\Pi^0_1$ class in which all members satisfy a property. A non-zero Turing degree ${\bf a}$ is {\em minimal} if there is no degree ${\bf b}$ such that ${\bf 0<b<a}$. Relativizing the minimal degree construction to any set, we know that every degree {\em has} a minimal cover, yet not every degree {\em is} a minimal cover. Our motivation arises from the work of Groszek and Slaman \cite{GS97} and from the problem that which Turing degrees {\em are} minimal covers, particularly their relationship with members of $\Pi^0_1$ classes. We consider a question asked by Andy Lewis-Pye in a personal communication, whether or not there exists a non-empty $\Pi^0_1$ class in which no member is a minimal cover, and we give a positive answer, thus showing that degrees of minimal covers cannot be a basis for $\Pi^0_1$ classes.



\section{Notation and Terminology}

We shall first give our notation and then give some background knowledge on $\Pi^0_1$ classes. We assume some familiarity with basic properties of relative computability and Turing degrees. For a detailed account of computability, the reader may refer to \cite{Soare},\cite{Cooper}, or \cite{Downey}. Readers who are familiar with computability theoretic notions and $\Pi^0_1$ classes may skip to Section \ref{sec:overall}.

Let $\omega$ denote the set of natural numbers. We let $2^{<\omega}$ denote the set of all finite sequences of 0's and 1's. We denote sets of natural numbers by uppsercase Latin letters $A,B,C$. 

The subset relation (not necessarily proper) is denoted by $\subset$. We identify a set $A\subset\omega$ with its characteristic function $f:\omega\rightarrow\{0,1\}$ such that, for any $n\in\omega$, if $n\in A$ then $f(n)=1$; otherwise $f(n)=0$. We let $\{\Psi_i\}_{i\in\omega}$ be an effective enumeration of the Turing functionals. Turing functionals will also be denoted by uppercase Greek letters $\Psi,\Phi,\Theta,\Xi$, etc. We say that $\Psi_e$ is {\em total} if it is defined on every argument, otherwise it is called {\em partial}. The {\em join} of any given two sets $A$ and $B$ is denoted by $A\oplus B=\{2i: i\in A\}\cup\{2i+1: i\in B\}$. $\Psi_e(A;n)\downarrow=m$ denotes that the $e$-th Turing functional with oracle $A$ on argument $n$ is defined and equal to $m$. $\Psi_e(A;n)\uparrow$ denotes it is not the case that $\Psi_e(A;n)\downarrow$. We will write $\Psi_e(A;n)[s]$ to mean $\Psi_e(A;n)$ as defined at stage $s$. $A\leq_T B$ means $B$ computes $A$ via some Turing functional $\Psi_e$. So $A\leq_T B$ iff there exists some $e\in\omega$ such that $A=\Psi_e(B)$. $A<_T B$ means $A\leq_T B$ and $B\not\leq_T A$. If $A\leq_T B$ and $B\leq_T A$, then $A$ an $B$ are {\em Turing equivalent}, which is denoted by $A\equiv_T B$. The {\em Turing degree} of a set $A\subset\omega$ is the set $\{X : X\equiv_T A\}$. We denote degrees with boldcase letters ${\bf a,b,c}$. 

We denote finite strings in $2^{<\omega}$ by lowercase Greek letters $\sigma,\tau,\eta,\rho,\pi,\upsilon$. They may also be denoted by the same letters with accents and subscripts such as $\sigma',\sigma'',\sigma^-,\sigma_0,\sigma_1$, etc. We let $\sigma*\tau$ denote the concatenation of $\sigma$ followed by $\tau$. For two strings $\sigma$ and $\tau$, we let $\sigma\subset\tau$ denote that $\sigma$ is an initial segment of $\tau$ as a substring. We say a string $\sigma$ is {\em compatible} with $\tau$ if either $\sigma\subset\tau$ or $\tau\subset\sigma$. Otherwise we say that $\sigma$ and $\tau$ are {\em incompatible}. Similarly, we say that $\sigma$ {\em extends $\tau$} if $\tau\subset\sigma$, i.e., $\sigma$ is an {\em extension} of $\tau$. Let $|\sigma|$ denote the length of $\sigma$. The unique string of length 0 is called the {\em empty string} and it is denoted by $\lambda$. We let $\sigma(i)$ denote the $(i+1)$st bit of $\sigma$.

For any $\sigma\in 2^{<\omega}$ and $n\in\omega$, conventionally we let $\Psi_e(\sigma;n)$ be defined and equal to $\Psi_e(A;n)$ if $\sigma(i)=A(i)$ for all $i<|\sigma|$ and if computing $\Psi_e(A;n)$ requires only values $A(i)$ for $i<|\sigma|$. Let $A\upharpoonright z$ and $\sigma\upharpoonright z$ denote, respectively, the restriction of $A(x)$ or $\sigma(x)$ to those $x\leq z$. 
For a set $A\subset\omega$, we define the {\em jump} of $A$, denoted by $A'$, to be the set $\{e:\Psi_e(A;e)\downarrow\}$. We write the $n$-th jump of a degree ${\bf a}$ as ${\bf a}^{(n)}$

A set $T$ of strings is {\em downward closed} if $\sigma\in T$ and $\tau\subset\sigma$ implies $\tau\in T$. Occasionally we refer to downward closed sets of strings as {\em trees}. We say that a set $A$ {\em lies on} a tree $T$ if there exist infinitely many $\sigma$ in $T$ such that $\sigma\subset A$. A set $A$ is a {\em path} on $T$ if $A$ lies on $T$. We denote the set of infinite paths of $T$ by $[T]$.  
We say that a string $\sigma\in T$ is {\em infinitely extendible} if there exists some $A\supset\sigma$ such that $A\in[T]$. 
If $\sigma,\tau\in T$ and $\sigma\subset\tau$ and there does not exist $\sigma'$ with $\sigma\subset\sigma'\subset\tau$ then we say that $\tau$ is an {\em immediate successor} of $\sigma$ in $T$ and $\sigma$ is the {\em immediate predecessor} of $\tau$ in $T$. 

We say that $\mathcal{P}\subset 2^\omega$ is a {\em $\Pi^0_1$ class} if there exists a downward closed computable set of strings $T$ such that $\mathcal{P}=[T]$. We can then have an effective enumeration $\{\Lambda_i\}_{i\in\omega}$ of downward closed computable sets of strings such that for any $\Pi^0_1$ class $\mathcal{P}$ there exists some $i\in\omega$ such that $\mathcal{P}$ is the set of all infinite paths through $\Lambda_i$. 

\subsection{Background on $\Pi^0_1$ classes}

In this section, we will overview some works in the literature to give an overall idea of how our claim is related to other works. However, for a detailed survey on $\Pi^0_1$ classes we refer the reader to \cite{CenzerRec} and \cite{Diamondstone}. A $\Pi^0_1$ class is an effectively closed subset of $2^\omega$ Cantor space. One important property of $\Pi^0_1$ classes is that for any recursively axiomatizable theory (the deductive closure of a recursively enumerable set of sentences in a language), the set of complete consistent extensions can be seen as a $\Pi^0_1$ class \cite{ShoenfieldAxiom}. The opposite direction is also proved in \cite{Ehrenfeucht}. That is, any $\Pi^0_1$ class can be seen as the set of complete consistent extensions of an axiomatizable theory. The compactness property of the Cantor space is provided by Weak K\"{o}nig's Lemma which tells us that if $T$ is an infinite downward closed set of finite strings, then there exists an infinite path through $T$. 

Countable $\Pi^0_1$ classes are another type of effectively closed sets. It is important to note that all countable $\Pi^0_1$ classes contain an isolated point and that every isolated point is computable \cite{Kreisel}. So if a $\Pi^0_1$ class contains no computable member then it must be uncountable and so $T$ must be {\em perfect}, i.e., every $\sigma\in T$ has at least two incompatible extensions in $T$.

Effectively closed sets in which no two members are of the same degree are called {\em $\Pi^0_1$ choice classes}, and the properties of these classes were studied in \cite{Cevik2016}. Any countably infinite $\Pi^0_1$ class has members of the same degree, hence $\Pi^0_1$ choice classes cannot be countably infinite. They also cannot be finite unless it is a singleton since all members of a finite $\Pi^0_1$ must be computable. Another property of $\Pi^0_1$ choice classes is that they do not contain any 1-random set. 

We are particularly interested in complexity of members of $\Pi^0_1$ classes in the Turing degree universe. Some of the most important and frequently used results are {\em basis} theorems: a basis theorem tells us that every non-empty $\Pi^0_1$ class has a member of a particular kind. Anything which is not a basis is called {\em non-basis}. The most celebrated {\em Low Basis Theorem} of Jockusch and Soare \cite{JS1} states that every non-empty $\Pi^0_1$ class contains a member of low degree, i.e., a degree ${\bf a}$ such that ${\bf a'=0'}$. Same authors proved that any non-empty $\Pi^0_1$ class contains a member of hyperimmune-free degree, i.e., a degree ${\bf a}$ such that for any $A\in{\bf a}$ and for any function $f\leq_T A$, there exists a computable function $g$ such that $g(n)\geq f(n)$ for all $n$. These results were proved by forcing with $\Pi^0_1$ classes in which we successively move from a set to one of its subsets in order to {\em force} satisfaction of a given requirement. Another important basis theorem for $\Pi^0_1$ classes is that every non-empty $\Pi^0_1$ class has a member of recursively enumerable (r.e.) degree, in particular, the leftmost path of any downward closed computable set of strings is of r.e. degree. One interesting result by Jockusch and Soare \cite{JS1} is that every $\Pi^0_1$ class which does not contain a recursive member contains members of degrees ${\bf a}$ and ${\bf b}$ such that ${\bf a\wedge b=0}$. However, this does not hold for the cupping case. In fact, it was shown in \cite{Cevik2021} that there exists a non-empty $\Pi^0_1$ class $\mathcal{P}$ with no recursive member such that $\emptyset'\not\leq_T A\oplus B$ for any $A,B\in\mathcal{P}$. Another non-basis result, given in \cite{JS2}, is that the class of r.e. degrees strictly below ${\bf 0'}$ does not form a basis. Similarly, the class of recursive sets does not form a basis either since there exists a $\Pi^0_1$ class such that all members are non-recursive. We call $\Pi^0_1$ classes with no recursive member {\em special} $\Pi^0_1$ classes. In \cite{Diamondstone}, it was proven that every non-empty special $\Pi^0_1$ class contains a member of properly low$_n$ degree, i.e., a degree ${\bf a}$ such that ${\bf a}^{(n)}={\bf 0}^{(n)}$ but ${\bf a}^{(n-1)}\neq{\bf 0}^{(n-1)}$. We say that a degree ${\bf a}$ satisfies the {\em join property} if for all non-zero ${\bf b<a}$ there exists ${\bf c<a}$ such that ${\bf b\vee c=a}$. In \cite{Cevik2021}, it was shown that there exists a non-empty special $\Pi^0_1$ class in which no member satisfies the join property. 

An {\em antibasis} theorem tells us that a $\Pi^0_1$ class cannot have all/any members of a particular kind without having a member of every Turing degree. Kent and Lewis \cite{KL} proved the {\em Low Antibasis Theorem} which says that if a $\Pi^0_1$ class contains a member of every low degree, then it contains a member of every degree. In \cite{Cevik}, the latter result was extended in relation to the jump hierarchy, that for a given degree ${\bf a\geq 0'}$, if a $\Pi^0_1$ class $\mathcal{P}$ contains members of every degree ${\bf b}$ such that ${\bf b'=a}$, then $\mathcal{P}$ contains members of every degree. A local version of this result is also given in the same work. That is, when ${\bf a}$ is also $\Sigma^0_2$, it suffices in the hypothesis to have a member of every $\Delta^0_2$ degree ${\bf b}$ such that ${\bf b'=a}$.

\section{Overall idea of the proof}\label{sec:overall}

The question of whether there exists a non-empty special $\Pi^0_1$ class in which no member is a minimal cover was asked by Andy Lewis-Pye in a personal communication around 2013. Our aim is to show that degrees of minimal covers cannot be a basis for $\Pi^0_1$ classes. There are three closely related results in the literature. First was provided by Lewis \cite{LewisPi01StrongMinimalHFP}, where he showed that there exists a non-empty special $\Pi^0_1$ class every member of which is of degree with strong minimal cover. Second related result was given by Sasso \cite{Sasso1970}, where he modified Sacks' minimal degree construction below ${\bf 0'}$ to obtain, recursively in ${\bf 0'}$, a tree the limit of every branch of which is recursive or of minimal degree. Of course Sasso didn't really construct a $\Pi^0_1$ class, but it can be thought of as a limit computable version of an effectively closed set in which all members satisfy the minimal degree requirements. The third relevant work is a paper by Groszek and Slaman \cite{GS97}, where they show the existence of a non-empty special $\Pi^0_1$ class in which all members compute a set which is of minimal degree. In fact, their theorem uses some level of non-uniformity, for that they construct a $\Pi^0_1$ class in which every member is either of minimal degree or computes a set of minimal degree. The non-uniformity here is used to satisfy the second literal of the disjunction. They define a co-r.e. tree $M$, with no terminal nodes, such that the $n$-th minimality requirement
\[
\Psi_n(A)\leq_T A\Rightarrow \Psi_n(A)\textrm{ is recursive or }A\leq_T \Psi_n(A)
\]
is satisfied for all $A\in[M]$ except for countably many members $B\in[M]$. They start with Sacks' \cite{Sacks1961} minimal degree construction method below ${\bf 0'}$, that is, by enumerating a $\Psi_n$-splitting tree $T_n$ in $M$ and guaranteeing that all branches through $M$ are either branches through $T_n$, in which case we have $A\leq_T\Psi_n(A)$, or they extend a terminal node of $T_n$. Now for branches $B$ extending a terminal node of $T_n$, they guarantee that $B$ computes a non-recursive r.e. set via some reduction that the authors construct. Hence by Yates' theorem \cite{Yates1970}, that every non-recursive r.e. set computes a set of minimal degree, $B$ computes a set which is of minimal degree. Of course, by the hyperimmune-free basis theorem such a $\Pi^0_1$ class must contain a member of minimal degree since hyperimmune-free sets cannot compute any non-zero r.e. set. Thus, the $\Pi^0_1$ class constructed by Groszek and Slaman must necessarily contain a member of minimal degree.

For our problem, the most naive approach in constructing a non-empty special $\Pi^0_1$ class with the property that no member is a minimal cover would be to make sure that all members are non-recursive and of r.e. degree. But this cannot work since every non-empty $\Pi^0_1$ class contains a member of hyperimmune-free degree. Of course if we drop the non-recursiveness condition, then obviously a $\Pi^0_1$ class with all members recursive gives an effectively closed set in which no member is a minimal cover. But this does not give solution to our problem. The approach we will consider to solve our problem will be similar to the Groszek-Slaman construction, but will contain some different features. So first we shall discuss to what extent they are related and how they differ from each other so as to give a general idea of our strategy. Like Groszek and Slaman, we will find it more convenient to build a $\Pi^0_1$ class from its complement, i.e., we define a co-r.e. tree $T$ with no terminal nodes such that no member of $[T]$ is a minimal cover. Groszek and Slaman use the {\em level $n$ construction} for satisfying the $n$-th minimality requirement as the main procedure of their overall construction with a subroutine what they call the {\em $\Psi_n$-subconstruction} procedure, where they take care of the paths $B$ for which they fail to satisfy the condition $B\leq_T\Psi_n(B)$ for countably many sets $B$ in the $\Pi^0_1$ class, but instead they guarantee that $B$ computes a non-recursive r.e. set via a Turing reduction that they also construct (unless $\Psi_n(B)$ is partial or recursive, thus $B$ is subject to a successful level $n+1$ construction). In other words, when enumerating a $\Psi_n$-splitting tree, as long as the outcome is $\Pi_2$, i.e., in the $\Psi_n$-splitting region of their tree at stage $s$, the $n$-th minimality requirement is satisfied by ensuring $B\leq_T\Psi_n(B)$. In the $\Sigma_2$ outcome, where no useable $\Psi_n$-splitting are found until stage $s$, they ensure that either $B$ computes a non-recursive r.e. set, or $\Psi_n(B)$ is partial or recursive and $B$ is subject to a successful level $n+1$ minimality construction. 

\begin{figure}[ht]
\begin{center}
\includegraphics{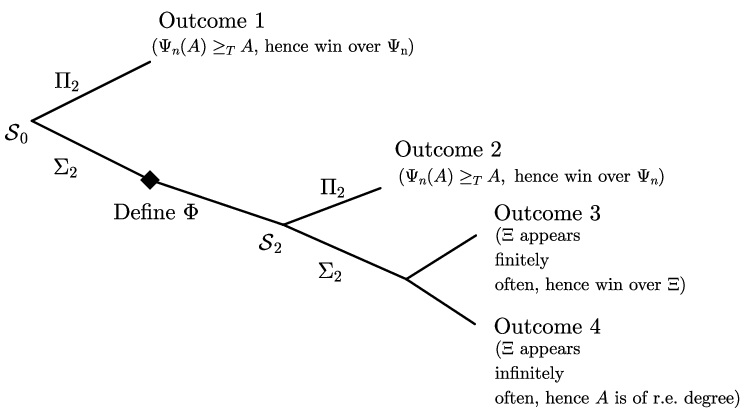}
\end{center}
\caption{Outcomes and wins.}
\label{fig:figureoutcomes}
\end{figure}

We will define $\Psi_n$-splittings inside our tree at each stage, just like in Groszek and Slaman's construction, and our approach for the $\Pi_2$ outcome of the minimality requirement (in the $\Psi_n$-splitting region of the tree) is the same. If $A$ is a branch on a $\Psi_n$-splitting tree $T$, then $A \leq_T \Psi_n(A)$. The proof of this is quite standard. Given $\Psi_n(A)$ we can generate increasingly long segments of $A$ recursively in $\Psi_n(A)$. Given $\sigma\subset A$, assuming that $A$ lies on the tree $T$, either $\sigma*0$ or $\sigma*1$ must be an initial segment of $A$, and we have to decide which one of them is. Since $T$ is $\Psi_n$-splitting, there exists some $x\in\omega$ such that $\Psi_n(\sigma*0; x)\downarrow\neq \Psi_n(\sigma*1; x)\downarrow$. But then only one of them can be compatible with $\Psi_n(A;x)$. We just take the one which is compatible, and this determines which of the two strings is an initial segment of $A$. So if $A \leq_T \Psi_n(A)$, the theorem is satisfied rather trivially and we will win by that particular $\Psi_n$. As for the $\Sigma_2$ outcome, though, our strategy will be different since our problem is not to ensure that the branches compute a set of minimal degree, but to ensure that if $A\in [T]$, where $T$, say, is the co-r.e. tree to be constructed, then there exists a reduction $\Phi$ satisfying that $\Psi(A)<_T \Phi(A)<_T A$. In other words, we carry out a density argument inside the region of $T$ where and as long as we fail to find any $\Psi_n$-splittings. This is where most of work in the paper will be given. The $n$-th minimality requirement, taken together with the non-recursiveness requirement that all branches through the tree are non-recursive, will roughly correspond to Groszek and Slaman's level $n$ construction. We will not work with the $n$-th minimality requirement aiming to satisfy it on all branches of $T$, as that would be an attempt to produce a $\Pi^0_1$ class all of whose members are of minimal degree, which is impossible anyway since every non-empty $\Pi^0_1$ class contains a member of r.e. degree (and r.e. degrees cannot be minimal due to Sacks' density theorem \cite{SacksDensity}). Instead, if we can't seem to get $A\leq_T\Psi_n(A)$ for satisfying the theorem trivially, we will succeed in other ways, such as making relevant functionals in the requirements (listed below) partial or we will make $A$ of r.e. degree (see Figure \ref{fig:figureoutcomes} for a sketch of possible outcomes).

Along the construction, we define $\Psi_n$-splittings inside the tree $T$ at each stage, and the $\Psi_n$-splitting region will grow over the dummy region where we have not yet found any $\Psi_n$-splittings for computations up to that stage. The density strategies (will be denoted by $\mathcal{S}_1$ and $\mathcal{S}_2$) will be divided into $\Pi_2$ and $\Sigma_2$ cases, depending if the node equipped with the strategy is in the $\Psi_n$-splitting region or otherwise of $T$ of that particular stage. In the $\Psi_n$-splitting region, we will have no problem in satisfying the theorem in a trivial way. The density requirements will hold trivially in the $\Pi_2$ case since we get $A\leq_T\Psi_n(A)$ for every $A$ that lies on the $\Psi_n$-splitting subtree of $T$. If, however, we observe that there are no $\Psi_n$-splittings by stage $s$, that is where we need the density requirements. In that case we start to implement a density construction for this $\Sigma_2$ region until we get $\Psi_n$-splittings, and we have to ensure that if $A$ is a branch that leaves a terminal node of the $\Psi_n$-splitting tree, then we guarantee that $\Phi(A)$ is total if $\Psi_n(A)$ is total. This is to ensure that $\Psi_n(A)\leq_T\Phi(A)$. We will observe that if it were the case that $\Theta(\Psi(A))=\Phi(A)$, then either $\Theta$ is partial (and so we win by that particular $\Theta$) or we can get an anti-chain of $\Phi$-splits on a unique infinite path $A$ where $\Psi$ is partial everywhere but on $A$. In that case we will show $A$ is of r.e. degree. 
Similarly, as for the other density requirement, if it were the case that $\Xi(\Phi(A))=A$ for some Turing functional $\Xi$, then we make sure that either $\Xi$ is partial (and so we win by that particular $\Xi$) or $A$ is of r.e. degree. 

\begin{definition}
Given any stage $s$, let $\sigma\in T_s$, where $T_s$ denotes the set of strings that are not enumerated out of $T$ before stage $s$. For a given $n\in\omega$, we say that $\sigma$ has an {\em active $\Psi_n$-split} at stage $s$ if there exist two incompatible strings $\sigma_1,\sigma_2\supset\sigma$ in $T_s$ such that $|\sigma_1 | = |\sigma_2 | = s$ and $\Psi_n(\sigma_1)\neq\Psi_n(\sigma_2)$.
\end{definition}

\noindent The requirements will be as follows. For any $A\in[T]$,
\vspace{0.5cm}

\noindent $\mathcal{S}_0(n):$ Either $\Psi_n(A)\geq_T A$ or $\exists\sigma\subset A$ $\exists s\in\omega$ such that $\sigma$ has no active $\Psi_n$-split extension after stage $s$.

\noindent $\mathcal{S}_1(n)$: $\Theta(\Psi_n(A))=\Phi(A)\Rightarrow$ $A$ is of r.e. degree.

\noindent $\mathcal{S}_2(n)$: $\Xi_n(\Phi(A))=A\Rightarrow$ $A$ is of r.e. degree.

\noindent $\mathcal{D}(n)$: $A\neq\Psi_n(\emptyset)$.
\vspace{0.5cm}

The non-recursiveness requirements $\mathcal{D}$ will be handled inside the $\mathcal{S}_1$ and $\mathcal{S}_2$ strategies separately. The strategy to satisfy the requirement denoted by $\mathcal{S}_0$, when taken together with the strategy to satisfy the $\mathcal{D}$ requirements, will act as our ``level $n$ construction" akin to the one in Groszek-Slaman construction. The {\em density requirements}, where we act so long as we don't find any active $\Psi_n$-splittings, are denoted by $\mathcal{S}_1$ and $\mathcal{S}_2$, and will have wins by $\Pi_2$ or $\Sigma_2$ outcomes. The density requirements make sure there is a set inbetween $\Psi_n(A)$ and $A$ for $n\in\omega$. We sketch the possible outcomes in Figure \ref{fig:figureoutcomes} and we will go through it in the next section. The outcomes for $\mathcal{S}_1$ (not shown in Figure \ref{fig:figureoutcomes}) will be similar, except that instead of $\Xi$, we will have $\Theta$ partial. How we satisfy the $\mathcal{S}_1$ and $\mathcal{S}_2$ requirements will become clear once we define the strategies.

In the next section we give the proof of our main result based on the idea we discussed so far. It should be noted that the construction mechanics is complicated enough that it may not be possible to work out with a single level of procedure and single layer of indices for construction stages, but we may require to work with multiple construction layers, as generally occurs in recursive constructions alike.

\section{Proof of the main theorem}\label{sec:proof}

We shall now give the main result of this paper.

\begin{definition}
A set $A$ is a {\em minimal cover} for a set $B<_T A$ if there is no $C$ such that $B<_T C<_T A$. A degree ${\bf a}$ is a {\em minimal cover} for a degree ${\bf b<a}$ if there is no ${\bf c}$ such that ${\bf b<c<a}$.
\end{definition}

\begin{theorem}\label{thm:main}
There exists a non-empty special $\Pi^0_1$ class in which no member is a minimal cover for any set. 
\end{theorem}

The proof is a type of density argument inside a tree. We construct a non-empty special $\Pi^0_1$ class $\mathcal{P}$ such that for any $\Psi(A)$ of any given $A\in\mathcal{P}$, if $\Psi(A)<_T A$, we define $\Phi(A)$ with $\Phi(A)\leq_T A$ and $\Psi(A)\leq_T \Phi(A)$, and we satisfy the requirements given above in Section \ref{sec:overall}.
\vspace{0.5cm}

Our approach in satisfying $\mathcal{S}_0$ (together with $\mathcal{D}$ requirements) is similar to the method discussed in Groszek and Slaman. One difference however is the application of the density requirements used in the splitting trees and that the $\mathcal{D}$ requirements are handled separately by $\mathcal{S}_1$ and $\mathcal{S}_2$ strategies.

We shall represent the $\Pi^0_1$ class as the collection of branches through a co-r.e. tree $T$ with no terminal strings. For this we describe a recursive construction that will determine $T$ by specifying at each stage $s$, which strings of the full binary tree will be enumerated out of $T$ before stage $s$. We denote by $T_s$ the set of strings that are not enumerated out of $T$ before stage $s$ and we let $T=\bigcap_s T_s$. We ensure that each $T_s$ is an infinite perfect tree, i.e., every string in $T_s$ has at least two incompatible extensions. In the end, $T$ will be an infinite tree with no terminal strings such that $[T]\neq\emptyset$.
\vspace{0.5cm}

\noindent {\bf Notation.} For two strings $\sigma$ and $\tau$ in $T$ such that $\sigma\subset\tau$, let us write $\sigma\sim\tau$ to denote that there is no string $\eta\in T$ such that $\sigma\subset\eta$ and that $\eta$ and $\tau$ are incompatible. For convenience, when we enumerate a string $\tau$ out of $T$, we also enumerate out any $\sigma\in T$ such that either $\tau\subset\sigma$ or $\sigma\subset\tau$ and $\sigma\sim\tau$. This preserves that $T_s$ has no terminal strings. Enumerating out of $T$ finitely many strings (and their extensions) and guaranteeing to leave at least one string will ensure all together that each $T_s$ is an infinite perfect tree.

\subsection{Control mechanics}

We give some more conventions regarding the construction. When working with trees, the terms `string' and `node' will often be used interchangably to denote the elements of the tree. The overall recursive construction will contain the instructions of how to define the $\Pi^0_1$ class $\mathcal{P}=[T]$ for a co-r.e. tree $T$. But the main procedure used in the overall construction is the {\em level $n$ construction}, i.e., the $\mathcal{S}_0(n)$-strategy. The level $n$ construction inside $\Psi_n$-splitting tree above the node $\tau$ describes a part of the overall recursive construction which produces above $\tau$ a sequence of trees $\{T_s\}_{s\in\omega}$. The $\Psi_n$-splitting tree will be a splitting tree on $T$. The level $n$ construction is what corresponds to the strategy for $\mathcal{S}_0$ (when taken together with the strategy for non-recursiveness $\mathcal{D}$ requirements). The density strategies will act as a subprocedure in the level $n$ construction so long as we don't find active $\Psi_n$-splits. For the level $n$ construction, our notation for handling the construction mechanics will be slightly different from that of the Groszek-Slaman level $n$ construction in the sense we will explain below. 


During the construction we will place on nodes various strategies for the requirements that we listed earlier. 
A node may be equipped with more than one strategy. A fixed priority ordering can be defined on the equipped strategies to execute them in some particular order. 
At any stage $s$, each node which is equipped with a strategy {\em acts} when the hypothesis of the associated requirement holds up to stage $s$. 

We also use what may be called {\em flagging} along the construction in the following sense. At each stage $s$, every node will have a binary {\em flag} value, namely $\Pi$ and $\Sigma$, associated with the $n$-th Turing functional $\Psi_n$ for each $n$ such that $n\leq s$. Since the outcomes of the $\mathcal{S}_0$ requirement can be divided into $\Sigma_2$ and $\Pi_2$ cases, the flag value of the node denotes whether the node is in the splitting region or otherwise of $T_s$ with respect to $\Psi_n$ at stage $s$. We say that a node $\sigma$ is in the {\em $\Pi$-region} in $T_s$ with respect to $\Psi_n$ if for any $\sigma_1,\sigma_2\supset\sigma$ of length $\leq |s|$ in $T_s$ we have that $\Psi_n(\sigma_1;x)[s]\neq\Psi_n(\sigma_2;x)[s]$ for some $x\in\omega$. If there is no active $\Psi_n$-splittings above $\sigma$, then $\sigma$ is in the {\em $\Sigma$-region} of $T_s$. Of course a node may be in different regions with respect to different indices of Turing functionals. For example if every incompatible extensions of $\sigma$ are $\Psi_n$-splitting but there are no active $\Psi_m$-splittings above $\sigma$, then $\sigma$ is in the $\Pi$-region with respect to $\Psi_n$, but in the $\Sigma$-region with respect to $\Psi_m$. Depending on which region $\sigma$ is in, with respect to an index $i$, the relevant versions of the $\mathcal{S}_1$ and $\mathcal{S}_2$ strategies will be executed accordingly. The $\Pi$ and $\Sigma$ versions of the $\mathcal{S}_i$-strategy, for $i\in\{1,2\}$, are denoted respectively as $\mathcal{S}_{i,\Pi}(n)$ and $\mathcal{S}_{i,\Sigma}(n)$. In the given scenario, then, $\mathcal{S}_{1,\Pi}(n)$ and $\mathcal{S}_{2,\Pi}(n)$ will be executed above $\sigma$ in $T_s$ along with $\mathcal{S}_{1,\Sigma}(m)$ and $\mathcal{S}_{2,\Sigma}(m)$. 

The control mechanics of the strategy for $\mathcal{S}_2$ may need a bit more care. We will use a subroutine in the $\mathcal{S}_2$ strategy which enumerates and adjusts possible $\Phi$ values in such a way to satisfy some condition related to $\Psi$ values that we will explain shortly. 

In constructions alike, it is customary to describe which requirements, in particular $\mathcal{S}_2$ requirements in this case, are `active at a node $\sigma$' at any given stage of the construction and which ones deserve to be called `complete' in order to move on to the next requirement. 

Suppose that $\sigma$ is a node and $i\in\omega$ is a given index. Let $\sigma_0$, $\sigma_1$ be two incompatible successors of $\sigma$ on which some strategies are placed. We describe the control mechanics taking $\mathcal{S}_2$ as an example, but it can be similarly applied to the others. We say that an $\mathcal{S}_2(i)$-strategy is {\em active at $\sigma$} if $\sigma$ is not yet enumerated out of $T$ by stage $s$ and $\Xi(\Phi(A))[s]=\tau\subset A$ for some functional $\Xi$. When the $\mathcal{S}_2(i)$-strategy is active at $\sigma$, given $\Xi_i$, the $\mathcal{S}_2(i)$-strategy placed on $\sigma$ will search for each $k\in\{0,1\}$,
\vspace{0.5cm}

an extension $\sigma^*_k\supset\sigma_k$ of length $\leq |s|$ in $T_s$ such that $\Xi_i(\Phi(\sigma^*_k))\supset\sigma_k$.
\vspace{0.5cm}

\noindent When such $\sigma^*_k$ is found, we say that $\sigma$ is {\em complete} (for the associated strategy) for all pairs $(i,\sigma')$ such that $\sigma'\supset\sigma_k$. 
A node with $\mathcal{S}_2(i)$-strategy being complete means that the hypothesis of the $\mathcal{S}_2(i)$ requirement is now satisfied for $\Psi_i$ and extensions of $\sigma_k$ at stage $s$. Once a node is complete for a strategy, we can pass the control to the one with next lower priority. 

The {\em level} of a string $\sigma$ in $T_s$ will refer to the number of proper initial segments of $\sigma$ in $T_s$. 
We decide whether or not $i\in\omega$ requires the attention of $\mathcal{S}_2$ at node $\sigma$ as follows. Let us say that $i$ requires {\em $\mathcal{S}_2$-attention} at $\sigma$ unless there exists some proper initial segment $\rho\subset\sigma$ with an $\mathcal{S}_2(i)$-strategy such that $\mathcal{S}_2(i)$ is active at $\rho$ but $\rho$ is not complete for all $(i,\tau)$ such that $\tau\supset\rho$. Let us also denote the set of $\mathcal{S}_2$ requirements active at a node $\sigma$ by $\beta_\sigma$. This will be determined by their indices. For example if $\beta_\sigma=\{m\}$, then only $\mathcal{S}_2(m)$ is active at $\sigma$. 

If $\sigma=\lambda$, then $\beta_\sigma=\{0\}$.

Suppose that $\sigma\neq\lambda$ and $\sigma$ is a string of level $n>0$. Let $\sigma'$ be the initial segment of $\sigma$ equipped with an $\mathcal{S}_2$-strategy placed and which is of level $n-1$. If $\sigma'$ is $\mathcal{S}_2$-complete for all $(i,\sigma)$ such that $i\in \beta_{\sigma'}$, then

\begin{center}
$\beta_\sigma=\beta_{\sigma'}\cup\{i'\}$,
\end{center}

\noindent where $i'$ is the least such number not in $\beta_{\sigma'}$ which requires $\mathcal{S}_2$-attention at $\sigma$. Otherwise, $\beta_\sigma$ is the set of all $i\in \beta_{\sigma'}$ which require $\mathcal{S}_2$-attention at $\sigma$. 



\vspace{0.5cm}

Strategy for $\mathcal{S}_1$ is executed over a pair of indices $(i,j)$, for arbitrary functionals $\Theta$ and $\Psi$. However, without loss of generality, we may assume for simplicity working with a single index $n=\langle i,j\rangle$ for either of the arbitrary functionals appearing in $\mathcal{S}_1$ requirements. Whenever an $\mathcal{S}_1(n)$-strategy is placed on a string $\sigma$ such that $\sigma_0$ and $\sigma_1$ are two incompatible extensions of $\sigma$, it will search for strings $\sigma^*_0\supset\sigma_0$ and $\sigma^*_1\supset\sigma_1$ satisfying either 

\begin{center}
$\sigma_0=\Phi(\sigma_0^*)\subset\Theta(\Psi_n(\sigma_0^*))$\hspace{1cm} or\hspace{1cm} $\sigma_1=\Phi(\sigma_1^*)\subset\Theta(\Psi_n(\sigma_1^*))$.
\end{center}



For $\mathcal{D}$ requirements, given $n\in\omega$, the $\mathcal{D}(n)$-strategy will be executed to satisfy the non-recursiveness requirement above $\sigma$. That is, it will ensure that for any $A\supset \sigma$ in $[T]$, $A\neq\Psi_n(\emptyset)$.
\vspace{0.5cm}

There will be finitely many requirements active on a given node bounded by stage $s$. At each stage, the strategy placed on a node performs the instructions for all of these in order of some fixed priority. 
We consider the $\mathcal{S}_0(n)$-strategy acting as the main procedure of the overall construction of our $\Pi^0_1$ class since $\mathcal{S}_1$ and $\mathcal{S}_2$ strategies will have two versions, each depending on the outcome of $\mathcal{S}_0$ at each stage. Similarly, the $\mathcal{D}(n)$-strategy will be performed in accordance with $\mathcal{S}_1$ and $\mathcal{S}_2$.

We also define a set $\Phi^+$ of strings for maintaining possible $\Phi$ axioms during the construction. We define axioms of the form $\Phi(\sigma)\supset\tau$. Set of possible axioms of $\Phi$ that are enumerated in $\Phi^+$ may not be necessarily consistent with each other when they are considered as axioms, but we will be able enumerate in the construction proper axioms of $\Phi$ from $\Phi^+$. 
The important point when enumerating axioms for $\Phi$ is that we will aim to make sure that $\Phi(A)$ is total if $\Psi(A)$ is total for any $A\in [T]$.

Now we give the strategies. First we describe how the strategy for $\mathcal{S}_0$ works above a given string $\sigma$.

\subsection{$\mathcal{S}_0(n)$-strategy above $\sigma$}   

\begin{figure}[ht]
\begin{center}
\includegraphics{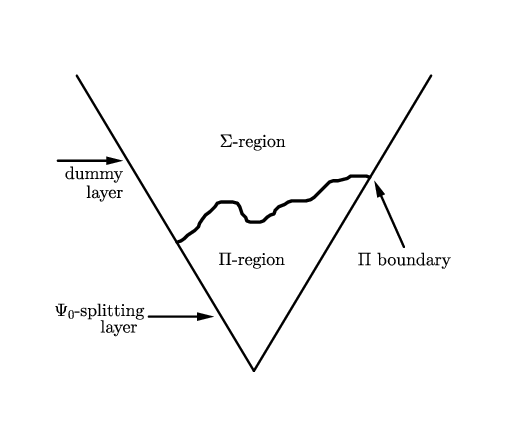}
\end{center}
\caption{The $\Pi$-region grows over the $\Sigma$-region.}
\label{fig:figureregions}
\end{figure}

Given stage $s$ and index $n\in\omega$, we say that a string $\sigma\in T$ is a {\em $\Pi_{s,n}$-boundary point} if $\sigma$ is in the $\Pi$-region of $T_s$ but every $\tau\supset\sigma$ is in the $\Sigma$-region of $T_s$ with respect to $\Psi_n$. Intuitively, $\Pi_{s,n}$-boundary points determine the longest strings in $T$, as defined at stage $s$, above which there are no active $\Psi_n$-splittings (see Figure \ref{fig:figureregions} for depiction).

The role of $\mathcal{S}_0$ will be to update the $\Pi$-region of $T$ at each stage. The strategy will define the $\Pi_{s,n}$-boundary points at each stage $s$ which will determine the layers containing $\Psi_n$-splitting strings and strings above which there is no active $\Psi_n$-split. When we say $\sigma$ is a $\Pi_{s,n}$-boundary point, we mean this so relative to the given index $n$ in the strategy and the stage $s$ we are in. So occasionally we will drop the indices and just write $\Pi$-boundary point as it should be understood relativized to that stage and functional index. Same thing when we say a string is in the ``$\Pi$-region", where it should be understood that, in fact, we mean this so relative to a functional index $n$.

Let $\sigma$ be a $\Pi_{s,n}$-boundary point on which $\mathcal{S}_0(n)$-strategy is placed. Intuitively, at stage $s$, in $\mathcal{S}_0(n)$-strategy we define $T_s\subset 2^{<\omega}$ such that either $\Psi_n(\sigma_1)[s]\neq\Psi_n(\sigma_2)[s]$ for every incompatible extensions $\sigma_1\supset\sigma$ and $\sigma_2\supset\sigma$ of length $s$ in $T_s$ or that there exists some $\tau\supset\sigma$ above which there are no active $\Psi_n$-splittings of length $s$. 
\vspace{0.5cm}

\noindent Instructions of step $s$ inside $T_s$ above $\sigma$ are as follows:


At stage $0$: Define $\tau\supset\sigma$ as the least string in $T$ such that $\tau*0$ and $\tau*1$ are not yet enumerated out.


At stage $s>0$, suppose that $\sigma$ is a $\Pi_{s,n}$ boundary point on which an $\mathcal{S}_0(n)$-strategy is placed, where $n\leq s$. We see if there exist strings $\sigma_1,\sigma_2\supset\sigma$ of length $s$ which are extendible in $T_s$ such that $\Psi_n(\sigma_1)[s]\neq\Psi_n(\sigma_2)[s]$. If not, we do nothing with the $\Pi$-boundary points and just place $\mathcal{S}_1(n)$ and $\mathcal{S}_2(n)$-strategies on $\sigma$ if not placed yet. Otherwise, when we find $\sigma_1$ and $\sigma_2$, we enumerate out all $\tau\in T$, for $\tau\supset\sigma$, such that $\tau$ is incompatible with $\sigma_1$ and $\sigma_2$. 
Take two least extensions $\tau_{1,i}$ and $\tau_{2,i}$ of $\sigma_i$, for each $i\in\{1,2\}$, and declare them as the new $\Pi$-boundary points. Then place $\mathcal{S}_1(n)$ and $\mathcal{S}_2(n)$-strategies on two incompatible extensions of $\tau_{1,i}$ and $\tau_{2,i}$.

\vspace{0.5cm}

Since the $\Pi$-boundary points of the previous stage are extended by the ones of the current stage, this ensures that the $\Pi$-region expands over the $\Sigma$-region. We also ensure that $\mathcal{S}_0(n)$ does not terminate any nodes in the $\Pi$-region or the permanent $\Sigma$-region. This can be shown by induction on $s$ that for every $n$ and $\sigma$, step $s$ of level $n$ construction inside $T_s$ above $\sigma$ does not enumerate $\sigma$ out of $T$. At step $s+1$ of the construction, we execute step $t$ for some $t\leq s$ of level $n+1$ construction inside $T_t$ above $\tau\supset\sigma$, and we apply the inductive hypothesis.

By the end of stage $s$, $\mathcal{S}_0(n)$ will have two outcomes in $T_s$. Either we find $\Psi_n$-splittings up to length $s$ above the given string $\sigma$ in the hypothesis, or there is no active $\Psi_n$-splittings. Let us denote these outcomes, respectively, as $\Pi_2$ and $\Sigma_2$.

$\mathcal{S}_1$ and $\mathcal{S}_2$ will act accordingly for each outcome of $\mathcal{S}_0$. Then, the strategy $\mathcal{S}_1$ has two versions, namely $\mathcal{S}_{1,\Pi}$ and $\mathcal{S}_{1,\Sigma}$. Similarly for $\mathcal{S}_2$. Suppose that $\sigma$ is a node with an $\mathcal{S}_0(i)$-strategy and consider the $\Pi_2$ outcome of $\mathcal{S}_0(i)$. In that case, we carry out the construction as if $\sigma$ has an active $\Psi_i$-split. So all low priority $\mathcal{S}_1$ and $\mathcal{S}_2$ requirements will perform their $S_{2,\Pi}(i)$-strategy. For instance, whenever we have a $\Pi_2$ outcome of $\mathcal{S}_0(i)$ and $\Sigma_2$ outcome of $\mathcal{S}_0(i+1)$, all low priority $\mathcal{S}_1$ and $\mathcal{S}_2$ requirements will be performing in that region their $S_{j,\Pi}(i)$ and $S_{j\Sigma}(i+1)$ strategies, for $j\in\{1,2\}$, meaning that $\mathcal{S}_1$ and $\mathcal{S}_2$ will act as if there are active $\Psi_i$-splits but no active $\Psi_{i+1}$-splits above the string on which the strategy is placed.

$\mathcal{S}_1$ and $\mathcal{S}_2$ strategies will therefore either end up with a $\Pi_2$ or $\Sigma_2$ outcome depending on the flag value of the equipped node. In each case we have a different win. In Figure \ref{fig:figureoutcomes} we sketch these possible outcomes and wins.

Note that the $\mathcal{S}_0(i)$-strategy we gave above is for a fixed Turing functional $\Psi_i$. For other functionals, without loss of generality we may use, just as in the Sacks construction, $\Psi_j$-splittings inside a $\Psi_i$-splitting tree for $j>i$. In the overall construction, as we will give in the end, the $\mathcal{S}_0(i+1)$-strategy will enumerate $\Psi_{i+1}$-splittings inside the $\Psi_i$-splitting subtree of $T$. 


Let us now analyze Figure \ref{fig:figureoutcomes}. Outcome 1 is achieved when we suppose that we are in the $\Psi_n$-splitting region. In this case we have $\Psi_n(A)\geq_T A$ for $A$ on the $\Psi_n$-splitting subtree of $T$, so we have a win over that particular $\Psi_n$. Otherwise, we have subcases. Let us suppose that we are in the $\Sigma$-region produced by the $\mathcal{S}_0(n)$-strategy. We suppose there exists some $\sigma\subset A$ such that there is no active $\Psi_n$-splittings above $\sigma$. This is where we need to define $\Phi$ when $\Psi_n$ is defined. Assume that we are working with the $\mathcal{S}_2$-strategy. Now $\mathcal{S}_2$ will both work in the $\Pi$-region and $\Sigma$-region of $T_s$. In the $\Pi$-region, we have the Outcome 1 (=Outcome 2), where we get $\Psi_n(A)\geq_T A$, and we have a win over the functional $\Psi_n$. Outcome 3 is obtained when $\mathcal{S}_2$ has a $\Sigma_2$ outcome and when $\Xi$ appears only finitely many times. In this case we have a win over $\Xi$ for the fact that $\Xi$ being partial. Otherwise, we have Outcome 4, in which case we have that $\Psi$ is partial on all paths except a unique infinite path $A$ which we can argue is of r.e. degree, hence cannot be a minimal cover. Similar condition holds for $\mathcal{S}_1$, only that we replace $\Xi$ with $\Theta$.
\vspace{0.5cm}

\noindent Before we give the strategy for $\mathcal{S}_1$, let us first explain the strategy for $\mathcal{S}_2$. The strategies for $\mathcal{S}_1$ and $\mathcal{S}_2$ are invoked inside the (main) $\mathcal{S}_0$-strategy. The $\mathcal{S}_1$ and $\mathcal{S}_2$ strategies will execute their instructions depending if the node equipped with that strategy is in the $\Pi$ or $\Sigma$-region of $T_s$. For the $\mathcal{S}_2$ strategy, for instance, if the flag value of $\sigma$ is $\Pi$, we execute $\mathcal{S}_{2,\Pi}(n)$; otherwise we execute $\mathcal{S}_{2,\Sigma}(n)$.

\subsection{$\mathcal{S}_2(n)$-strategy above $\sigma$}     


We always try to diagonalize if possible. If we find out that $\Xi(\Phi(\sigma))\neq\sigma$ at any stage of any version of the $\mathcal{S}_2$ strategy given below, we enumerate out from $T$ every node incompatible with $\sigma$ and cancel all strategies placed on them.

\subsubsection{$\mathcal{S}_{2,\Pi}(n)$-strategy instructions above $\sigma$} 

If the flag value of $\sigma$ is $\Pi$, that is, if we are working in a $\Psi_n$-splitting region above $\sigma$, then $\Psi_n(\sigma)$ computes $\sigma$ so things work out easier and the requirements are satisfied trivially. We just need to proceed with the non-recursiveness requirement to make sure no infinite path above $\sigma$ is recursive. We place the $\mathcal{D}(n)$ strategy on two incompatible strings of $\sigma$ and execute $\mathcal{D}(n)$ for both. At most one of them will get diagonalized and enumerated out of $T$. Then we place $\mathcal{S}_2(n+1)$ strategy on the other string which remains in $T$. This ensures that the $\mathcal{S}_{2,\Pi}(n)$-strategy does not enumerate out $\sigma$ from $T$.

\subsubsection{$\mathcal{S}_{2,\Sigma}(n)$-strategy instructions above $\sigma$} 

Suppose that we are at stage $s$ and the flag value of $\sigma$ with respect to $\Psi_n$ is $\Sigma$, and so there are no active $\Psi_n$-splits above $\sigma$. Let $\tau_0$ and $\tau_1$ be two incompatible extensions of $\sigma$ that are not yet enumerated out of $T$. We ensure the following property:
\vspace{0.5cm}

\noindent For every $\tau_0'\supset\tau_0$ and $\tau_1'\supset\tau_1$, if $\Psi_n(\tau'_0)$ and $\Psi_n(\tau'_1)$ are defined and compatible with each other, then $\Phi(\tau_0')$ and $\Phi(\tau_1')$ are compatible.\hfill $(\star)$
\vspace{0.5cm}

Let $\sigma_0\supset\sigma$ be a node on which the next highest priority $\mathcal{D}(n)$-strategy is placed. Let $\sigma_1$ and $\sigma_2$ be two successors of $\sigma_0$. We wait until stage $s'$ such that $\Psi_n$ is defined on all strings $\tau$ in $T_{s'}$ up to length $s'$, where $\tau\supset\sigma_i$, for $i\in\{1,2\}$. Until we find so, the idea is to place $\Phi$-splits on every string which is a $\Pi$-boundary point of stage $s$ (This will allow us to make $\Psi$ partial on branches where we have $\Phi$-splittings). We can do this so as follows. Suppose $\eta\supseteq\sigma$ is a $\Pi$-boundary point at stage $s$ with respect to functional $\Psi_n$. See if there exist two incompatible strings $\eta_0$, $\eta_1$ in $\Phi^+$ extending $\Phi(\eta)$ whenever $\Phi(\eta)$ is defined (if $\Phi(\eta)$ is undefined at this stage, we do nothing). If they exist, choose $\eta_0'$ and $\eta_1'$ extending $\eta$ such that $\eta_0'$ and $\eta_1'$ are incompatible. We then place $\mathcal{S}_2(n+1)$-strategies on $\eta_0'$ and $\eta_1'$ and enumerate axioms $\Phi(\eta_0')\supset\eta_0$ and $\Phi(\eta_1')\supset\eta_1$ for a choice of shortest possible length of strings $\eta_0$ and $\eta_1$. If they do not exist, let $\tau'\supset\Phi(\eta)$ be the longest extension of $\Phi(\eta)$ in $\Phi^+$. Choose two incompatible strings $\eta_0$ and $\eta_1$ extending $\eta$. Place $\mathcal{S}_2(n+1)$-strategies on $\eta_0$ and $\eta_1$, and enumerate the axioms $\Phi(\eta_0)\supset\tau*0$ and $\Phi(\eta_1)\supset\tau*1$. Also enumerate $\tau'*0$ and $\tau'*1$ into $\Phi^+$. 
\vspace{0.5cm}

Along the construction we may run into the problem of having either one of the extensions of $\sigma$ removed from $T$, which may interfere with the preservation of $(\star)$. To avoid this problem and to preserve $(\star)$, we use `reflections' in the sense that we arrange $\Phi$ values in such a way that any $\Phi$ value above $\sigma_1$ will be available in $\Phi^+$ above $\sigma_2$. For this, we modify a subroutine which appeared in an earlier work \cite{Cevik2021}. The $\Phi^+$-regulation subprocedure given below will ensure that if the branch extending one of the incompatible extensions of $\sigma$, say $\eta$, is enumerated out from $T$ at a later stage, we will still be able to use the other $\Psi_n$-splitting pair via $\Phi^+$ to define the value of $\Phi(\eta)$. So after enumerating $\tau'*0$ and $\tau'*1$ into $\Phi^+$, we run the $\Phi^+$-regulation procedure on $\sigma$.
\vspace{0.5cm}

\noindent{\bf Subroutine for $\Phi^+$-regulation procedure on argument $\sigma$}: 

The idea is to preserve compatible $\Phi$ pairs above incompatible extensions of $\sigma$. Let $\sigma_1,\sigma_2\supset\sigma$ be the immediate incompatible successors of $\sigma$. We ensure that if $\tau\supset\sigma_1$, then there exists some $\pi\supset\sigma_2$ such that either $\Phi(\tau)\subset\Phi(\pi)$ or $\Phi(\pi)\subset\Phi(\tau)$ if $\Phi(\tau)$ and $\Phi(\pi)$ are defined in $\Phi^+$.

Pick the longest string $\eta\in\Phi^+$, where $\sigma\subset\eta$, such that $\Phi(\eta)$ is defined and all $\eta_0\subset\eta$ are complete for all $i\in\beta_{\eta_0}$. For every $\rho$ of length $\leq|\eta|$ and incompatible with $\eta$, we define an anti-chain of strings in $\Phi^+$ whose $\Phi$ values are compatible with $\Phi(\eta)$. For this we take $\Phi(\eta)$. Then for every $\rho$ of length $l\leq |\eta|$ such that $\rho\not\subset\eta$ and $\eta\not\subset\rho$, we enumerate $\rho$ into $\Phi^+$. Then we define the axiom $\Phi(\rho)\supset\Phi(\eta)*0^l$, where $0^l$ denotes $l$ consecutive $0's$ whenever $\Psi(\rho)$ and $\Psi(\eta)$ are defined and compatible. This ensures that all incompatible extensions $\rho,\eta\supset\sigma$ have compatible values in the $\Phi$ domain whenever $\Psi(\rho)$ and $\Psi(\eta)$ are compatible (More details as to why this preserves the $(\star)$ property are given in Lemma \ref{lem:total}). We then enumerate $\rho*0$ and $\rho*1$, as well as $\eta*0$ and $\eta*1$ into $\Phi^+$. This allows us to extend the possible domain of $\Phi$ for enumerating further axioms.

\noindent {\bf End of subroutine}
\vspace{0.5cm}

%
%
%
%

The $\Phi^+$-regulation subprocedure is to ensure that every value of $\Phi$ above $\sigma_1$ is also available in $\Phi^+$ as a value for $\Phi$ above $\sigma_2$ in case one of them gets diagonalized. We call these $\Phi$ values {\em reflections} of each other.

\begin{figure}[ht]
\begin{center}
\includegraphics[scale=0.70]{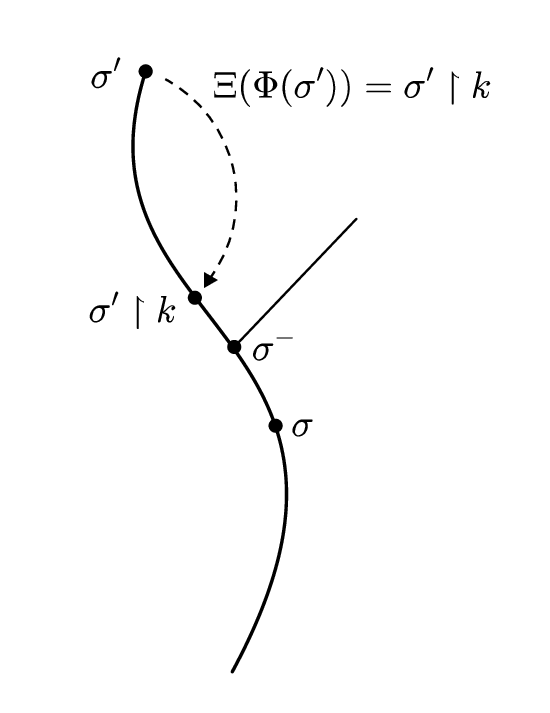}
\end{center}
\caption{Suppose that there are no active $\Psi_n$-splits above $\sigma$. Provided that $(\star)$ holds, if $\Psi_n(\sigma_0)$ is compatible with $\Psi_n(\sigma_1)$ for every $\sigma_0,\sigma_1\supset\sigma$ then $\Xi$ cannot be defined above any $\sigma^-$ incompatible with $\sigma'\upharpoonright k$.}
\label{fig:figureone}
\end{figure}

Now $\mathcal{S}_2(n)$ will search for a string $\sigma'_i\supset\sigma_i$ such that $\Xi(\Phi(\sigma'_i))\supset\sigma_i\upharpoonright k$ for some fixed $k\geq |\sigma |$. If there does not exist such string, we are fine as the requirement is satisfied trivially in this case. Otherwise, we declare every string above $\sigma_i\upharpoonright k$ and incompatible with $\sigma'_i$ to be terminal and enumerate them out from $T$. We also enumerate out those strings from $T$ whose $\Phi$ values are defined and are incompatible with $\Phi(\sigma'_i)$. We let $T_{s+1}$ to be the resulting tree. Note that if $\Xi(\Phi(\sigma'_i))\supset\sigma_i\upharpoonright k$ then $\Psi$ must be partial anywhere but above $\sigma'_i$ since otherwise $\Xi(\Phi^*(\sigma'_i))$ would compute $\sigma_{1-i}$ (see Figure 3). It is also worth noting that we ensure to keep at least one node alive above $\sigma$. The reason is that if a string $\sigma$ is removed for being inconsistent with the computation of $\Xi(\Phi)$ on a fixed argument, then we must be keeping strings incompatible with $\sigma$ since one of the two incompatible strings must be consistent with the latter computation.
\vspace{0.5cm}

Regardless of which version we follow, after these instructions we place $\mathcal{S}_2(n+1)$-strategies on strings of the least level on which no lower priority $\mathcal{S}_2$-strategy is placed yet.

\subsection{$\mathcal{S}_1(n)$-strategy above $\sigma$}   


We try to diagonalize if possible as usual. That is, whenever we find that $\Theta_n(\Psi_n(\sigma))\neq\Phi(\sigma)$, we declare every string incompatible with $\sigma$ to be terminal and enumerate them out from $T$, and so $\mathcal{S}_1(n)$ is satisfied above $\sigma$. Until we falsify the hypothesis of the $\mathcal{S}_1(n)$ requirement, the strategy is instructed to perform the following in accordance with the flag value of $\sigma$, i.e., depending if $\sigma$ is in the $\Sigma$ or the $\Pi$-region of $T_s$ with respect to $\Psi_n$.

\subsubsection{$\mathcal{S}_{1,\Pi}(n)$-strategy instructions above $\sigma$}  

Suppose that $\sigma$ has an active $\Psi_n$-split. Then we will have no problem in satisfying the requirement trivially since we automatically get $\Psi_n(A)\geq_T A$ and we have a win over that particular $\Psi$.

We nevertheless need to apply the non-recursiveness strategy to make sure no infinite path above $\sigma$ is recursive. We place the $\mathcal{D}(m)$-strategy on two incompatible strings $\tau_0, \tau_1$ such that  $\tau_0,\tau_1\supset\sigma$ and we execute $\mathcal{D}(m)$ on both for least such $m$ that has not been picked yet for the non-recursiveness requirement. At most one of them will get diagonalized, say $\tau_i$, and enumerated out of $T$. We then place $\mathcal{S}_2(n+1)$-strategy on $\tau_{1-i}$. Also define $\tau_{1-i}$ to be a new $\Pi$-boundary point.

\subsubsection{$\mathcal{S}_{1,\Sigma}(n)$-strategy instructions above $\sigma$}  

Suppose that an $\mathcal{S}_1(n)$-strategy is placed on $\sigma$ and that there are no active $\Psi_n$-splits above $\sigma$. Given $T_s$ at stage $s$ we perform the following instructions:

\begin{enumerate}
\item Fix some new witness $l\in\omega$.

\item See if there are two incompatible strings $\sigma_1$ and $\sigma_2$ in $T_s$, extending $\sigma$, such that $\Psi_n(\sigma_1)$ and $\Psi_n(\sigma_2)$ are defined up to $l$. If so, then

\begin{enumerate}
\item[(i)] Place a $\mathcal{D}(m)$-strategy on $\sigma_1$ and $\sigma_2$ for least $m$ that has not been picked yet, and define axioms for $\Phi$ such that $\Phi(\sigma_1)$ and $\Phi(\sigma_2)$ are incompatible up to $l$.

\item[(ii)] Keep both extensions $\sigma_1$ and $\sigma_2$ in $T_s$ until we see either one of two things happen:

\begin{enumerate}
\item[(a)] There exists an extension of $\sigma_i$, for some $i\in\{1,2\}$, say $\sigma'_i\supset\sigma_i$, such that $\Theta(\Psi_n(\sigma'_i))=\Phi(\sigma_i)$.
\item[(b)] The $\mathcal{D}(m)$-strategy decides to remove one of $\sigma_i$.
\end{enumerate}
\end{enumerate}
\end{enumerate}

\noindent Now one of the two cases may happen for $\Theta$ as the construction goes:
\vspace{0.5cm}

\noindent (i) $\Theta$ may appear finitely often.

\noindent (ii) $\Theta$ may appear infinitely often.
\vspace{0.5cm}

\noindent Case (i): Consider first the finite outcome. It may be that $\Theta$ stops appearing above $\sigma$. If we have a finite outcome, then we have a win over $\Theta$ and so there is nothing to prove as we can satisfy the $\mathcal{S}_1$ requirement above $\sigma$.
\vspace{0.5cm}

\begin{figure}[ht]
\begin{center}
\includegraphics[scale=0.80]{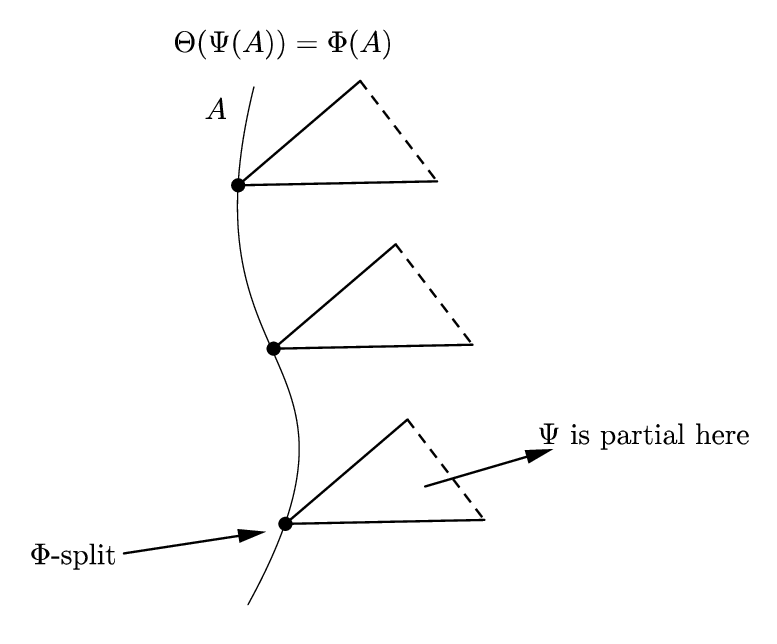}
\end{center}
\caption{If $\Psi_n$ is total on a unique infinite path $A$ and there are infinitely many $\Phi$-splits, $\Psi_n(A)\geq_T A$ and $\Phi(A)\geq_T\Psi_n(A)$.}
\label{fig:figuretwo}
\end{figure}

\noindent Case (ii): We follow a similar argument as in the $\Sigma_2$ outcome for $\mathcal{S}_2$. We keep both $\sigma_1$ and $\sigma_2$ extendible until the hypothsis in the $\mathcal{S}_1$ requirement is satisfied. If we find that $\Theta(\Psi_n(\sigma_i'))=\Phi(\sigma_i)$, we define all $\tau\supset\sigma_i$ to be terminal such that $\tau$ is incompatible with $\sigma_i'$. Note that when $\Theta$ of $\Psi_n$ is defined above for $\sigma$ which is compatible with $\sigma_1$, then $\Psi$ cannot be defined above for $\sigma_2$ because the range of $\Psi$ is below $\sigma_1$. If we later decide to kill the branch on which $\Theta$ is defined, then $\Theta$ may now get defined above $\sigma_2$. But let us suppose this does not happen, or at least it happens finitely many times. This means that $\Theta$ will be defined along a unique infinite path $A$ and $\Psi_n$ will be partial on every other path (see Figure \ref{fig:figuretwo}). So for this unique infinite path on which $\Theta$ is defined, we have $\Psi_n(A)\geq_T A$ and $\Phi(A)\geq_T A$. This infinite path, we will prove in Lemma \ref{lem:re}, is of r.e. degree. 

\subsection{$\mathcal{D}(n)$-strategy above $\sigma$}

Finally we explain the strategy $\mathcal{D}(n)$ for the non-recursiveness requirements. $\mathcal{D}(n)$ works in accordance with $\mathcal{S}_1$ and $\mathcal{S}_2$. Hence, it will have two versions, one for $\mathcal{S}_1$ and the other for $\mathcal{S}_2$.
\vspace{0.5cm}


\noindent Instructions of $\mathcal{D}(n)$ for $\mathcal{S}_1$ at stage $s$:

Let $\tau$ be a string in $T_s$ equipped with a $\mathcal{D}(n)$-strategy and suppose that $\mathcal{D}(n)$ has not yet acted on $n$, in the sense that $n$ is a newly picked witness. We find the least $\tau'\supset\tau$ such that $\tau'$ is the least string on which no $\mathcal{S}_1(n)$-strategy is placed. Let $\tau'_0$ and $\tau'_1$ be two strings extending $\tau'$ such that $\tau'_0$ and $\tau'_1$ are incompatible. Define $\Phi(\tau'_0)$ and $\Phi(\tau'_1)$ to be two incompatible strings. If we find at some later stage that the $n$-th Turing functional $\Psi_n$  extends one of $\tau'_i$, for $i\in\{0,1\}$, then we declare all extensions of $\tau'_i$ which are incompatible with $\tau'_{1-i}$ as terminal, enumerate them out of $T$ and remove all strategies from these strings. Then place $\mathcal{D}(n+1)$-strategy on two incompatible extensions of $\tau'_{1-i}$.

\vspace{0.5cm}

\noindent Instructions of $\mathcal{D}(n)$ for $\mathcal{S}_2$ at stage $s$:

Suppose that a $\mathcal{D}(n)$-strategy is placed on $\sigma$. We want the $\mathcal{D}(n)$-strategy to decide which one of the two incompatible extensions of $\sigma$ to keep in $T$ at stage $s$. When working with $\mathcal{S}_2$ requirements we have to be careful that the domain of $\Phi$ should not exceed the domain of $\Psi_n$ before we decide which path to choose for diagonalization. Let $\sigma_0$ and $\sigma_1$ be two incompatible extensions of $\sigma$ in $T_s$. In order to keep at least one of $\sigma_0$ and $\sigma_1$ extendible, we leave reflections of $\Phi$ values above both $\sigma_0$ and $\sigma_1$ in case the $\mathcal{D}(n)$-strategy at a later stage decides to remove one of them and we make sure that $\Phi$ is not defined until $\Psi_n$ gets defined. For this it suffices to ensure that it is not the case that every extendible $\sigma'_0\supset\sigma_0$ has a $\Phi$ value incompatible with its reflection. More precisely, we make sure that above $\sigma_0$, there exist some $\sigma'_0$ for which there exists some $\sigma'_1\supset\sigma_1$ such that $\Phi(\sigma'_0)$ and $\Phi(\sigma'_1)$ are compatible if $\Psi_n(\sigma_0')$ and $\Psi_n(\sigma_1')$ are compatible. This is ensured by the $(\star)$ property, which is preserved by the $\Phi^+$-regulation procedure. So we only enumerate $\Phi$ axioms while leaving a reflection of either of $\Phi(\sigma_0')$ or $\Phi(\sigma_1')$ in $\Phi^+$. Again, we do this because we want to preserve $(\star)$ for making $\Phi$ values of extensions compatible whenever the $\Psi$ values are compatible. If it were the case that all $\Phi(\sigma'_0)$ were incompatible with $\Phi(\sigma'_1)$, then assuming that $\Phi(\sigma'_0)\downarrow$ and $\Phi(\sigma'_1)\downarrow$, and that there are no active $\Psi_n$-splits above $\sigma$, we would not be able to make the reflections compatible with each other which is necessary for the preservation of $(\star)$. The rest of the instructions is similar to the instructions of $\mathcal{D}(n)$ for $\mathcal{S}_1$.


\subsection{Main construction}
We are now ready to define the main recursive construction, which is just initiated by enumerating a $\Psi_0$-splitting subtree of $2^{<\omega}$. 
\vspace{0.5cm}

{\bf Stage $0$.} Define $T_0=2^{<\omega}$. We define $\lambda$ to be the $\Pi$-boundary point initially. We let $\Phi^+=\{\lambda\}$ and enumerate $\Phi(\lambda)\supset\lambda$.


{\bf Stage $s+1$.} For any $\sigma$ of length $s+1$ in $T_{s+1}$, carry out step $s$ of the level $0$ construction inside $T_0$ above $\sigma$.

The construction generates, at each stage, a $\Psi_{n+1}$-splitting subtree over the $\Psi_n$-splitting path, on the regions where $\Psi_n$ is partial.

We let $T=\bigcap_s T_s$ be the co-r.e. tree produced by the construction. Now, the construction will never enumerate $\lambda$ out of $T$. Furthermore, since no terminal nodes are left in $T_s$, $T$ must be infinite and so $[T]\neq\emptyset$ (see Lemma \ref{lem:nonempty}).


\section{Verification}
We first argue that the constructed $\Pi^0_1$ class is non-empty and contains no recursive members.

\begin{lemma}\label{lem:nonempty}
$[T]$ is non-empty and it does not contain a recursive member.
\end{lemma}
\begin{proof}
It should be clear from both versions of the $\mathcal{D}(n)$-strategy that if $A\in[T]$ and $A=\Psi_n$ for some $n$, then there is some stage $s$ where $\sigma\subset A$ gets enumerated out of $T_s$ when we witness $\sigma\subset\Psi_n(\emptyset)[s]$. Hence, if $A$ is recursive then $A\not\in[T]$. If $\sigma$ is a node with a $\mathcal{D}(n)$-strategy, the fact that we run the $\mathcal{D}(n)$-strategy on each incompatible extension $\sigma_0,\sigma_1\supset$ ensures that $\mathcal{D}(n)$ does not enumerate out $\sigma$ from $T$.

To prove that $[T]$ is non-empty, it is clear from the instructions that if $\sigma$ is a node equipped with an $\mathcal{S}_i$-strategy, for $i\in\{0,1,2\}$, then $\sigma$ is not enumerated out of $T$. Also when we enumerate out some $\tau$ from $T$, we also enumerate out of $M$ all $\sigma\in T$ such that either $\tau\subset\sigma$ or $\sigma\subset\tau$ and $\sigma\sim\tau$. The latter property preserves the property that $T_s$ has no terminal nodes. Then these two facts together ensure that $T_s$ is an infinite perfect tree. Hence, $T=\bigcap_s T_s$ is an infinite perfect tree, and so $[T]\neq\emptyset$ by compactness of Cantor space.

\end{proof}

Next we show that the domain of $\Psi$ does not exceed that of $\Phi$. This ensures that $\Phi(A)$ computes $\Psi_n(A)$ for any $n\in\omega$.

\begin{lemma}\label{lem:total}
For any $A\in [T]$, $\Phi(A)$ is total if $\Psi(A)$ is total.
\end{lemma}
\begin{proof}
Assume that $A\in [T]$ and $\Psi(A)$ is total. Suppose that there exists some $\sigma$ and a stage $s$ after which $\Phi(\tau)$ is undefined for all $\tau\supset\sigma$ at all later stages $s'\geq s$. 
The totality of $\Phi(A)$ is ensured as follows whenever we define an axiom for $\Phi$: We take the least $n\in\omega$ such that $\Psi(\sigma;n)\uparrow$. We then make sure that for all $m<n$ and $j<m$, if $\Psi(\sigma;m)=\tau$ and if $\Psi(\sigma;j)=\eta$ such that $\eta\subset\tau$, then enumerate the axiom $\Phi(\sigma;n)=\tau$ for the longest string $\pi\in\Phi^+$ such that $\pi\subset\tau$. We then enumerate $\tau*0$ and $\tau*1$ into $\Phi^+$. If we follow this convention this runs contrary to what we supposed for a contradiction. Hence, if there exists some $k\in\omega$ for which $\Phi(\sigma;k)\uparrow$, then there exists some $l<k$ for which $\Psi(\sigma;l)\uparrow$.
\end{proof}

We now prove the last lemma needed to complete the proof of the main theorem. We shall first describe the success criteria for the execution of requirements. Let us say that, for each $i\in\{0,1,2\}$, {\em success conditions} of the $\mathcal{S}_i$ and $\mathcal{D}$ requirements depend on the satisfaction of two conditions:

\begin{enumerate}
\item[(i)] Step $0$ of the requirements are carried out. If step $s$ of the requirement is carried out, then there is a stage $t>s$ at which the requirement is carried out.
\item[(ii)] When a strategy is permanently satisfied after stage $t$, i.e., there is a stage $t$ such that the strategy is satisfied for every $t>s$, no string $\sigma$ above which the strategy is permanently satisfied is enumerated out of $T$.
\end{enumerate}

\begin{lemma}\label{lem:re}
Suppose that the success conditions for the $\mathcal{S}_i$-requirements, for each $i\in\{0,1,2\}$, are satisfied on infinite paths of $T$. Then any $A\in[T]$ in the $\Sigma$-region of $T$ for which $\Psi_n(A)$ is total is of r.e. degree.
\end{lemma}
\begin{proof}

The $\mathcal{S}_0(n)$-strategy enumerates a $\Psi_n$-splitting subtree of $T_s$. So with the $\mathcal{S}_0$-requirements we either end up with having a $\Pi_2$ outcome or $\Sigma_2$ outcome. Suppose that $A\in[T]$ and $A$ is in the $\Pi$-region of $T$, i.e., for every $\sigma\subset A$ there exist $\tau_1\supset\sigma$ and $\tau_2\supset\sigma$ in $T$ such that $\Psi_n(\tau_1;x)\downarrow\neq\Psi_n(\tau_2;x)$ for some $x\in\omega$. In this case, $\Psi_n(A)\geq_T A$ and so the theorem is satisfied trivially. We now prove that if $A$ lies in the permanent $\Sigma$-region of $T$ such that $\Psi_n(A)$ is total, then $A$ is of r.e. on that infinite branch of $T$. For the $\mathcal{S}_2$ requirements we have the following cases:

\begin{enumerate}
\item[(i)] If $\Theta$ appears finitely often, we have the Outcome 3 since in this case we win over $\Theta$.
\item[(ii)] If $\Theta$ appears infinitely often, we get Outcome 4. That is, we have that $\Psi_n(A)$ computes $A$ via $\Theta$ on a unique infinite path, and $\Psi_n$ is partial on any path incompatible with $A$. In this case, we argue below that $A$ is of r.e. degree.
\end{enumerate} 

We give the proof for the $\mathcal{S}_2$ requirements, but the same argument can be adapted to $\mathcal{S}_1$ requirements as well. Suppose that $A$ is the unique path in $[T]$ such that $\Psi_n(A)$ is total. We can construct an r.e. set $W$ which is Turing equivalent to $A$ as follows. Let $\mathcal{D}_\sigma$ denote the node $\sigma$ with a $\mathcal{D}(n)$-strategy placed on it. Define $W$ to be set of all $\mathcal{D}_\sigma$ such that there exists some stage $s$ where $\mathcal{D}_\sigma$ is in the $\Pi$-region with respect to $\Psi_n$ and $\mathcal{D}_\sigma$ is active at stage $s$. More precisely, we let
\[
W=\{\sigma\subset A:\exists s\text{ such that }\sigma\text{ is on the subtree with }\Psi_n(A)\geq_T A\text{ after stage } s\}.
\]



\noindent Clearly $W$ is an r.e. set since $W$ is a $\Sigma^0_1$ set. We claim that $W\equiv_T A$.

(i) The fact that $W\geq_T A$ is clear since $A$ is the unique infinite path in $W$.

(ii) We argue that $W\leq_T A$. We see if $\mathcal{D}_\sigma$ is used in the construction and so whether it remains in the tree at stage $|\sigma |$. If $\mathcal{D}_\sigma$ is used then our tree at any stage will look like an antichain of strings on a path $A$. Now suppose that $\sigma$ is a string from which we start our strategy. Given $\gamma\in 2^{<\omega}$, there are three cases to find out whether or not $\gamma\in W$. First case is that $\gamma$ could be off the whole tree. So if $\gamma$ is incompatible with every $\sigma$, then we know that $\gamma\not\in W$. Secondly, if it is the case that $\gamma\supset\sigma$ and $\gamma\subset A$, then $\gamma\in W$ and $\Psi(A)\geq_T A$ whenever $\gamma$ is on the infinite path on which $\Psi(A)$ is total. Finally, suppose that $\gamma\supset\sigma$ but $\gamma$ is incompatible with every $\tau\subset A$. In this case, we keep finding the node with the next highest priority $\mathcal{D}$-strategy until, at some stage $t$, we hit $\gamma$ and we see whether $\gamma$ is extended by an active node of length $t$. For convenience we may suppose that the extension has the same length as the length of the stage of the construction. So a string is in $W$ iff it is in $W[t]$ enumerated at stage $t$ where we put active node compatible with $A$. Then, $\mathcal{D}_\sigma\in W$ if and only if $\mathcal{D}_\sigma\in W$ by stage $|\sigma |$.
\end{proof}

This completes the proof of Theorem \ref{thm:main}. Hence, $[T]$ is a non-empty special $\Pi^0_1$ class in which no member is a minimal cover for any set. Thus, class of degrees of minimal cover is a not basis for $\Pi^0_1$ classes.

\bibliographystyle{amsplain}

\end{document}